\documentclass[12pt]{article}  
\def\sq{\hbox {\rlap{$\sqcap$}$\sqcup$}}
\def\sq{\hbox {\rlap{$\sqcap$}$\sqcup$}}
\def\R{ {\rm R \kern -.31cm I \kern .15cm}}
\def\C{ {\rm C \kern -.15cm \vrule width.5pt \kern .12cm}}
\def\Z{ {\rm Z \kern -.27cm \angle \kern .02cm}}
\def\N{ {\rm N \kern -.26cm \vrule width.4pt \kern .10cm}}
\def\1{{\rm 1\mskip-4.5mu l} }
\def\lsim{\raise0.3ex\hbox{$<$\kern-0.75em\raise-1.1ex\hbox{$\sim$}}}
\def\gsim{\raise0.3ex\hbox{$>$\kern-0.75em\raise-1.1ex\hbox{$\sim$}}}
\def\noi{\noindent}

\def\beq{\begin{equation}}   \def\eeq{\end{equation}}
\def\bea{\begin{eqnarray}}  \def\eea{\end{eqnarray}}

\def\noi{\noindent}
\def\beeq{\begin{eqnarray}} \def\eeeq{\end{eqnarray}}
\newcommand\mysection{\setcounter{equation}{0}\section}

\newcounter{hran}

\begin{document} 
\centerline{\large\bf Long range scattering for the Maxwell-Schr\"odinger system} 
 \vskip 3 truemm 
 \centerline{\large\bf  with arbitrarily large asymptotic data} 
  \vskip 0.8 truecm

\centerline{\bf J. Ginibre}
\centerline{Laboratoire de Physique Th\'eorique\footnote{Unit\'e Mixte de
Recherche (CNRS) UMR 8627}}  \centerline{Universit\'e de Paris XI, B\^atiment
210, F-91405 ORSAY Cedex, France}
\vskip 3 truemm

\centerline{\bf G. Velo}
\centerline{Dipartimento di Fisica, Universit\`a di Bologna}  \centerline{and INFN, Sezione di
Bologna, Italy}

\vskip 1 truecm

\begin{abstract}
We review the proof of existence and uniqueness of solutions of the Maxwell-Schr\"odinger system in a neighborhood of infinity in time, with prescribed asymptotic behaviour defined in terms of asymptotic data, without any restriction on the size of those data. That result is the basic step in the construction of modified wave operators for the Maxwell-Schr\"odinger system.
\end{abstract}

\vskip 1 truecm
\noi MS Classification :    Primary 35P25. Secondary 35B40, 35Q40.\par \vskip 2 truemm

\noi Key words : Long range scattering, Maxwell-Schr\"odinger system. \par 
\vskip 1 truecm

\noindent LPT Orsay 07-112\par
\noindent October 2007\par \vskip 3 truemm

\newpage
\pagestyle{plain}
\baselineskip 18pt
\mysection{Introduction}
\hspace*{\parindent} 
This paper is devoted to the theory of scattering and more precisely to the construction of modified wave operators for the Maxwell-Schr\"odinger (MS) system in space dimension 3, namely
\beq
\label{1.1e}
\left \{ \begin{array}{l} i\partial_t u = -(1/2) \Delta_A u + 
A_{e} u \\ \\ \sq A_{e} - \partial_t \left ( \partial_t A_{e} + \nabla \cdot A\right )  = |u|^2 \\ \\ \sq A + \nabla \left ( \partial_t A_{e} + \nabla \cdot A\right )  = {\rm Im}\ \overline{u} \nabla_A u\end{array} \right . \eeq

\noi where $u$ and $(A, A_e)$ are respectively a complex valued function and an ${I\hskip-1truemm R}^{3+1}$ valued function defined in space time ${I\hskip-1truemm R}^{3+1}$, $\nabla_A = \nabla - iA$ and $\Delta_A = \nabla_A^2$ are the covariant gradient and covariant Laplacian respectively, and $\sq = \partial_t^2 - \Delta$ is the d'Alembertian. An important property of that system is its gauge invariance, namely the invariance under the transformation
$$\left ( u, A, A_e \right ) \rightarrow \left (u \exp (- i \theta), A - \nabla \theta, A_e + \partial_t \theta \right ) \ ,$$

\noi where $\theta$ is an arbitrary real function defined in ${I\hskip-1truemm R}^{3+1}$. As a consequence of that invariance, the system (\ref{1.1e}) is underdetermined as an evolution system and has to be supplemented by an additional equation, called a gauge condition. Here,  we shall use exclusively the Coulomb gauge condition, $\nabla \cdot A = 0$. Under that condition, the equation for $A_e$ can be solved by
\beq
\label{1.2e}
A_e = - \Delta^{-1} |u|^2 = (4 \pi |x|)^{-1} \star |u|^2 \equiv g(u)
\eeq

\noi where $\star$ denotes the convolution in ${I\hskip-1truemm R}^3$. Substituting (\ref{1.2e}) and the gauge condition into (\ref{1.1e}) yields the formally equivalent system
$$\hskip 4 truecm \left\{ \begin{array}{ll}
i \partial_t u = - (1/2) \Delta_A u + g(u) u &\hskip 4.5 truecm (1.3)\\
&\\
\sq A = P \ {\rm Im} \ \overline{u} \nabla_A u &\hskip 4.5 truecm (1.4)
\end{array}\right .$$

\noi where $P = \1 - \nabla \Delta^{-1} \nabla$ is the projector on divergence free vector fields.\par

The MS system is known to be locally well posed both in the Coulomb gauge and in the Lorentz gauge $\partial_t A_e + \nabla \cdot A = 0$ in sufficiently regular spaces \cite{7r} \cite{8r}, to have weak global solutions in the energy space \cite{6r} and to be globally well posed in a space smaller than the energy space \cite{9r}.\par

A large amount of work has been devoted to the theory of scattering and more precisely to the existence of wave operators for nonlinear equations and systems centering on the Schr\"odinger equation, in particular for the Maxwell-Schr\"odinger system \cite{1r} \cite{2r} \cite{3r} \cite{5r} \cite{10r} \cite{11r}. As in the case of the linear Schr\"odinger equation, one must distinguish the short range case from the long range case. In the former case, ordinary wave operators are expected and in a number of cases proved to exist, describing solutions where the Schr\"odinger function behaves asymptotically like a solution of the free Schr\"odinger equation. In the latter case, ordinary wave operators do not exist and have to be replaced by modified wave operators including an additional phase in the asymptotic behaviour of the Schr\"odinger function. In that respect, the MS system in ${I\hskip-1truemm R}^{3+1}$ belongs to the borderline (Coulomb) long range case. General background and additional references on that matter can be found in \cite{3r} \cite{4r}.\par

The main step in the construction of the (modified) wave operators consists in solving the local Cauchy problem with infinite initial time and we shall concentrate on that step. In the long range case where that problem is singular, it amounts to constructing solutions with prescribed (singular) asymptotic behaviour in time. The results include the uniqueness and the existence of solutions with such prescribed asymptotic behaviour. The construction is performed by changing variables from $(u, A)$ to new variables, replacing the original system by an auxiliary system for the new variables, solving the corresponding problems for the auxiliary system and finally returning therefrom to the original one. The auxiliary system will be described in Section 2 below, and the uniqueness and existence results will be presented in Sections 3 and 4 respectively. The exposition is based on \cite{3r} \cite{5r}. \par

\mysection{The auxiliary system}
\hspace*{\parindent} 

In this section we derive the auxiliary system that will replace the original system (1.3) (1.4). We first replace the equation (1.4) by the associated integral equation with prescribed asymptotic data $(A_+, \dot{A}_+)$, namely
\beq
\label{2.1e}
A = A_0 - \int_t^{\infty} dt' \omega^{-1} \sin (\omega (t-t')) P \ {\rm Im} \left ( \overline{u} \nabla_A u \right ) (t')
\eeq

\noi where $\omega = (- \Delta )^{1/2}$, $A_0$ is the solution of the free wave equation $\sq A_0 = 0$ given by 
\beq
\label{2.2e}
A_0 = (\cos \omega t) A_+ + \omega^{-1} (\sin \omega t ) \dot{A}_+
\eeq

\noi and we assume $\nabla \cdot A_+ = \nabla \cdot \dot{A}_+ = 0$ in order to ensure the gauge condition $\nabla \cdot A_0 = 0$. We then perform a change of unknown functions which is well adapted to the study of the system (1.3) (\ref{2.1e}) for large time. The unitary group $U(t)$ which solves the free Schr\"odinger equation can be written as 
$$U(t) = \exp (i (t/2)\Delta ) = M(t) \ D(t)\ F\ M(t)$$

\noi where $M(t)$ is the operator of multiplication by the function
$$
M(t) = \exp \left ( ix^2/2t \right )\ ,
$$

\noi $F$ is the Fourier transform and $D(t)$ is the dilation operator defined by 
$$
D(t) = (it)^{-3/2} D_0(t) \qquad , \quad \left ( D_0 (t) f\right ) (x) = f(x/t)\ .
$$

\noi We first change $u$ to its pseudo conformal inverse $u_c$ defined by 
\beq
\label{2.3e}
u(t) = M(t)\ D(t) \ \overline{u_c (1/t)}
\eeq

\noi or equivalently
\beq
\label{2.4e}
\widetilde{u}(t) = \overline{F \widetilde{u}_c (1/t)}\ ,
\eeq

\noi where for any function $f$ of space time
\beq
\label{2.5e}
\widetilde{f}(t, \cdot ) = U(-t) \ f(t, \cdot )\ .
\eeq

\noi Correspondingly we change $A$ to $B$ defined by 
\beq
\label{2.6e}
B(t) = - t^{-1}D_0 (t)\ A(1/t)\ .
\eeq

\noi The transformation $(u, A) \to (u_c, B)$ is involutive. Furthermore it replaces the study of $(u, A)$ in a neighborhood of infinity in time by the study of $(u_c, B)$ in a neighborhood of $t =0$.\par

Substituting (\ref{2.3e}) (\ref{2.6e}) into (1.3) and commuting the Schr\"odinger operator with $MD$, we obtain
\begin{eqnarray*}
&&\left \{ \left (  i \partial_t + (1/2) \Delta_A - g(u) \right  ) u \right \} (t)\\
&&= t^{-2} \ M(t)\ D(t) \left \{ \overline{(i \partial_{t'} + (1/2) \Delta_{B(t')} - {\check B}(t') - t{'}^{-1} g(u_c(t')))u_c(t')} \right \}_{t'=1/t}
\end{eqnarray*}

\noi where for any ${I\hskip-1truemm R}^3$ vector valued function $f$ of space time
\beq
\label{2.7e}
{\check f}(t, x) = t^{-1} x \cdot f(t, x) \ .
\eeq

\noi Furthermore
$${\rm Im} \left ( \overline{u} \nabla_A u\right ) (t) = t^{-3} D_0(t) \left \{ x |u_c(t')|^2 - t' \left ({\rm Im}\ \overline{u}_c \nabla_B u_c\right ) (t') \right \}_{t'=1/t}$$

\noi by a direct computation, so that the system  (1.3) (\ref{2.1e}) becomes 
$$\hskip 2.1 truecm 
\left \{   \begin{array}{ll}
i \partial_t u_c = - (1/2) \Delta_B u_c + {\check B} u_c + t^{-1} g(u_c) u_c  &\hskip 4 truecm (2.8)\\
& \\
B_2 = {\cal B}_2 (u_c, B) &\hskip 4 truecm (2.9)
 \end{array} \right . 
$$

\noi where 
$$B_2 = B - B_0 - B_1 \ , \eqno(2.10)$$
$$B_0(t) = - t^{-1} D_0(t) \ A_0(1/t) \ , \eqno(2.6)_0$$
$$B_1 = B_1 (u_c) \equiv - F_1 (P x |u_c|^2)\ , \eqno(2.11)$$
$${\cal B}_2 (u_c, B) \equiv t\ F_2 \left ( P\ {\rm Im}\ \overline{u}_c \nabla_B u_c \right )\ , \eqno(2.12)$$
$$F_j(M) \equiv \int_1^{\infty} d\nu \ \nu^{-2-j}\ \omega^{-1} \sin ( \omega (\nu - 1))D_0(\nu) \ M(t/\nu)\ . \eqno(2.13)$$

\noi Here we take the point of view that $B_1$ is an explicit function of $u_c$ defined by (2.11) and that (2.10) is a change of dynamical variable from $B$ to $B_2$. The equation (2.9) then replaces (\ref{2.1e}).\par

In order to take into account the long range character of the MS system, we parametrize $u_c$ in terms of a complex amplitude $v$ and a real phase $\varphi$ by
$$u_c = v \exp (- i \varphi )\ . \eqno(2.14)$$

\noi The role of the phase is to cancel the long range terms in (2.8), namely the contribution of $B_1$ to ${\check B}$ and the term $t^{-1} g(u_c)$. Because of the limited regularity of $B_1$, it is convenient to split $B_1$ and $B$ into a short range and a long range part. Let $\chi \in {\cal C}^{\infty} ({I\hskip-1truemm R}^3, {I\hskip-1truemm R})$, $0 \leq \chi \leq 1$, $\chi (\xi ) = 1$ for $|\xi | \leq 1$, $\chi (\xi ) = 0$ for $|\xi | \geq 2$. We define
$$\hskip 2.1 truecm 
\left \{   \begin{array}{l}
{\check B}_L = {\check B}_{1L} = F^{\star} \chi (\cdot \ t^{1/2}) F\ {\check B}_1 \\
 \\
 {\check B}_S = {\check B}_0 + {\check B}_{1S} + {\check B}_2 \qquad , \quad {\check B}_{1S} = {\check B}_1 - {\check B}_{1L}\ ,
 \end{array} \right . \eqno(2.15)
$$

\noi We then obtain the following system for $(v, \varphi , B_2)$
$$\hskip 3.5 truecm 
\left \{   \begin{array}{ll}
i \partial_t v = Hv  &\hskip 5.5 truecm (2.16)\\
& \\
\partial_t  \varphi = t^{-1} g(v) + {\check B}_{1L} (v) &\hskip 5.5 truecm (2.17)\\
& \\
B_2 = {\cal B}_2 (v, K) &\hskip 5.5 truecm (2.18)
 \end{array} \right . 
$$

\noi where
$$H \equiv - (1/2) \Delta_K + {\check B}_{S}\ , \eqno(2.19)$$
$$K \equiv B + \nabla \varphi \ , \eqno(2.20)$$

\noi by imposing (2.17) as the equation for $\varphi$. Under (2.17), the equation (2.8) becomes (2.16). The system (2.16)-(2.18) is the auxiliary system which replaces the original system (1.3) (1.4).

\mysection{Uniqueness of solutions}
\hspace*{\parindent} 
In this section we present the uniqueness result at infinity in time for the MS system in the form (1.3) (\ref{2.1e}). Since the Cauchy problem for that system is singular at $t = \infty$, especially as regards the function $u$, the uniqueness result for that system takes a slightly unusual form. Roughly speaking it states that two solutions $(u_i, A_i)$, $i = 1,2$, coincide provided $u_i$ and $A_i - A_0$ do not blow up too fast and provided $u_1 - u_2$ tends to zero in a suitable sense as $t \to \infty$. In particular that result does not make any reference to the asymptotic data for $u$, which should characterize its behaviour at infinity.\par

We denote by $\parallel\ \cdot \ \parallel_r$ the norm in $L^r \equiv L^r ({I\hskip-1truemm R}^3)$, by $H^k$ the standard Sobolev spaces and by $\dot{H}^k$ their homogeneous versions, with $\dot{H}^1 \subset L^6$. We need the spaces 
\bea
\label{3.1e}
&&V = \left \{ v: v \in H^3,\  xv \in H^2 \right \}\ , \\
&&V_{\star} = \left \{ v:<x>^3 v \in L^2,\ <x>^2 \nabla v \in L^2 \right \}
\label{3.2e}
\eea

\noi which are Fourier transformed of each other, and  the dilation operator
\beq
\label{3.3e}
S = t \ \partial_t + x \cdot \nabla + 1 \ .
\eeq

\noi We denote nonnegative integers by $j$, $k$. \par

We first state the uniqueness result for the auxiliary system (2.16)-(2.18). \\

\noi {\bf Proposition 3.1.} {\it Let $0 < \tau \leq 1$, $I = (0, \tau ]$ and $\alpha \geq 0$. Let $A_0$ be a divergence free solution of the free wave equation $\sq A_0 = 0$ such that  $B_0$ defined by (2.6)$_0$ satisfy
\beq
\label{3.4e}
\parallel \nabla^k (t\partial_t)^j B_0 \parallel_{\infty}\ + \ \parallel \nabla^k {\check B}_0 \parallel_{\infty} \ \leq C\ t^{-k}
\eeq

\noi for $0 \leq j+k \leq 1$ and for all $t \in I$. Let $(v_i, \varphi_i, B_{2i})$, $i = 1,2$, be two solutions of the auxiliary system (2.16)-(2.18) such that $v_i \in L_{loc}^{\infty}(I, V)$, $\nabla \varphi_i \in L_{loc}^{\infty}(I, \dot{H}^1 \cap \dot{H}^2)$, $B_{2i } \in L_{loc}^{\infty} (I, \dot{H}^1)$, satisfying the estimates
\beq
\label{3.5e}
\parallel v_i (t) ; V \parallel\ \leq C(1 - \ell n \ t)^{\alpha}\ ,
\eeq
\beq
\label{3.6e}
\parallel \nabla B_{2i } (t) \parallel_2\ \leq C (1 - \ell n\ t)^{2\alpha}\ ,
\eeq
\beq
\label{3.7e}
\parallel <x> (v_1(t) - v_2(t))\parallel_2\ \leq C\ t^{1+\varepsilon} 
\eeq

\noi for some $\varepsilon > 0$ and for all $t \in I$, and such that $\varphi_1(t) - \varphi_2(t)$ tends to zero in $\dot{H}^1$ when $t\to 0$. Then $(v_1, \varphi_1, B_{21}) = (v_2 , \varphi_2, B_{22})$.}\\

\noi {\bf Remark.} Although each $\varphi_i$ separately is not in $\dot{H}^1$, it follows from the equation (2.17) and from the assumption (\ref{3.7e}) that $\varphi_1 - \varphi_2$ has a limit in $\dot{H}^1$ when $t\to 0$, which gives a meaning to the assumption that $\varphi_1 - \varphi_2$ tends to zero in $\dot{H}^1$ when $t \to 0$.\\

We next state the uniqueness result for the MS system in the form (1.3) (\ref{2.1e}). We recall that $\widetilde{u}(t) = U(-t) u(t)$. \\

\noi {\bf Proposition 3.2.} {\it Let $1 \leq T < \infty$, $I = [T, \infty )$ and $\alpha \geq 0$. Let $A_0$ be a divergence free solution of the free wave equation $\sq A_0 = 0$ satisfying
\beq
\label{3.8e}
\parallel \nabla^k S^j A_0(t) \parallel_{\infty}\ + \ \parallel \nabla^k x \cdot A_0(t) \parallel_{\infty} \ \leq C\ t^{-1} 
\eeq

\noi for $0 \leq j+k \leq 1$ and for all $t \in I$. Let $(u_i, A_i)$, $i = 1,2$, be two solutions of the system (1.3) (\ref{2.1e}) such that $\widetilde{u}_i \in L_{loc}^{\infty}(I, V_{\star})$, $A_i - A_0 \in L_{loc}^{\infty} (I, \dot{H}^1)$, satisfying the estimates
\beq
\label{3.9e}
\parallel \widetilde{u}_i (t) ; V_{\star} \parallel\ \leq C(1 + \ell n \ t)^{\alpha}\ , 
\eeq
\beq
\label{3.10e}
\parallel \nabla (A_i - A_0) (t) \parallel_2\ \leq C\ t^{-1/2} (1 + \ell n\ t)^{2\alpha}\ , 
\eeq
\beq
\label{3.11e}
\parallel <x/t> (u_1 - u_2)(t)\parallel_2\ \leq C\ t^{-1- \varepsilon} 
\eeq

\noi for some $\varepsilon > 0$ and for all $t \in I$. Then $(u_1, A_1) = (u_2 , A_2)$.}\\

The proof of Proposition 3.1 follows from an estimate of the difference of two solutions of the system (2.16)-(2.18), in particular of $v_1-v_2$ in $H^1 \cap F(H^1)$. Proposition 3.2 is proved by reducing it to Proposition 3.1. The main step of that reduction is the construction of the phases $\varphi_i$, and that step requires a strengthening of the assumptions in the form of the replacement of $\alpha$ by $7 \alpha + 3$ from (\ref{3.9e}) (\ref{3.10e}) to (\ref{3.5e}) (\ref{3.6e}). The assumptions (\ref{3.4e})  on $B_0$ and (\ref{3.8e}) on $A_0$ are equivalent. Since $SA_0$ and $x\cdot A_0$ are solutions of the free wave equations, those assumptions can be implemented under suitable well known conditions on the asymptotic data $(A_+, \dot{A}_+)$. Finally, Proposition 3.2 is sufficient to cover the case of the solutions of the MS system constructed in the next section.

\mysection{Existence of solutions}
\hspace*{\parindent} 
In this section, we sketch the proof of existence of solutions of the MS system (1.3) (\ref{2.1e}) with prescribed asymptotic behaviour at infinity in time, and for that purpose we first construct solutions of the auxiliary system (2.16)-(2.18) with prescribed asymptotic behaviour at time zero. Since that system is singular at $t = 0$, we proceed by an indirect method. We pick an asymptotic form for the dynamical variables, to be chosen later, we rewrite the auxiliary system in terms of the difference variables and we solve the latter for those variables under the condition that they tend to zero at time zero. Let $(v_a, \varphi_a , B_{1a}, B_{2a})$ be the asymptotic form of $(v, \varphi , B_{1}, B_{2})$ and define the difference variables
\beq
\label{4.1e}
\left ( w, \psi , G_1, G_2\right ) = \left ( v, \varphi, B_1, B_2 \right ) - \left (v_a, \varphi_a, B_{1a}, B_{2a} \right )\ ,
\eeq 

\noi where we do not assume that $B_{1a} = B_1(v_a)$ or $\varphi_a = \varphi (v_a)$, in order to allow for more flexibility. Substituting (\ref{4.1e}) into the auxiliary system (2.16)-(2.18), where in addition we reintroduce an equation for $B_1$, we obtain the following system for the difference variables~:
$$\hskip 2 truecm 
\left \{   \begin{array}{ll}
i \partial_t w = Hw + H_1 v_a - R_1  &\hskip 3 truecm (4.2)\\
& \\
\partial_t  \psi = t^{-1} g(w, 2v_a + w) + {\check B}_{1L} (w, 2v_a + w) -R_2 &\hskip 3 truecm (4.3)\\
& \\
G_1 = B_1 (w, 2v_a + w) - R_3&\hskip 3 truecm (4.4)\\
& \\
G_2 = {\cal B}_2 (w,2v_a+w, K) + b_2 (v_a, L) - R_4&\hskip 3 truecm (4.5)
 \end{array} \right . 
$$

\noi where $H$, $K$ are defined by (2.19) (2.20), $K_a = B_a + \nabla \varphi_a = B_0 + B_{1a} + B_{2a} + \nabla \varphi_a$, $L = G_1 + G_2 + \nabla \psi$, $B_1$ and ${\cal B}_2$ are defined by (2.11)-(2.13),
$$b_2(v_a, L) = {\cal B}_2 (v_a, K) - {\cal B}_2 (v_a, K_a)\ ,$$

\noi so that $b_2 (v_a, L)$ is linear in $L$,
$$H_1 = i L \cdot \nabla_{K_a} + (i/2) \Delta \psi + (1/2) L^2 + {\check G}_{1S} + {\check G}_{2}\ ,  \eqno(4.6)$$

\noi $ {\check G}_{1S}$ and  ${\check G}_{2}$ are defined according to (2.7) (2.15), and for any quadratic form $Q(f)$, we denote by $Q(f_1,f_2)$ the associated polarized sesquilinear form. The remainders $R_i$, $1 \leq i \leq 4$, are defined by
$$\hskip 1.5 truecm 
\left \{   \begin{array}{ll}
R_1 = i \partial_t v_a - H_a \ v_a=   i \partial_t v_a + (1/2) \Delta_{K_a} \ v_a - {\check B}_{aS}\ v_a &\hskip 3.2 truecm (4.7)\\
& \\
R_2 = \partial_t  \ \varphi_a -  t^{-1} g(v_a) - {\check B}_{1L} (v_a)  &\hskip 3.2 truecm (4.8)\\
& \\
R_3 = B_{1a}  - B_1(v_a) &\hskip 3.2 truecm (4.9)\\
& \\
R_4 = B_{2a} - {\cal B}_2 (v_a,K_a) \ . &\hskip 3 truecm (4.10)
 \end{array} \right . 
$$

\noi They measure the failure of the asymptotic form to satisfy the original system (2.16)-(2.18) and their time decay as $t \to 0$ measures its quality as an asymptotic form.\par

The resolution of the new system (4.2)-(4.5) now proceeds in two steps. \\

\noi \underline{Step 1.} One solves the system (4.2)-(4.5) for $( w, \psi , G_1, G_2)$ tending to zero as $t \to 0$ under general boundedness properties of $(v_a, \varphi_a , B_{1a}, B_{2a})$ and general decay assumptions on the remainders $R_i$, $1 \leq i \leq 4$, as $t \to 0$. For that purpose, one linearizes the system (4.2)-(4.5) partly into a system for new variables $(w', G'_2)$, namely
$$\hskip 1.5 truecm 
\left \{   \begin{array}{ll}
i \partial_t w' = Hw' + H_1 v_a - R_1  \\
& \\
G'_2 = {\cal B}_2 (w,2v_a+w, K) + b_2 (v_a, L) - R_4 \ .
 \end{array} \right . \eqno(4.11)
$$

\noi (There is no point in introducing new variables for $\psi$ and $G_1$ since they are explicitly given as functions of $w$ under the condition $\psi (0) = 0$). For given $(w, G_2)$ tending to zero as $t \to 0$, one solves the system (4.11) for $(w', G'_2)$ also tending to zero at $t \to 0$. This defines a map $\Gamma : (w, G_2) \to (w', G'_2)$. One then shows by a contraction method that the map $\Gamma$ has a fixed point. \\

\noi \underline{Step 2.} One constructs asymptotic functions $(v_a, \varphi_a , B_{1a}, B_{2a})$ satisfying the assumptions needed for Step 1. For that purpose, one solves the auxiliary system (2.16)-(2.18) approximately by an iteration procedure. It turns out that the second approximation is sufficient.\\

An essential ingredient in the implementation of Step 1 is the choice of the function space where the map $\Gamma$ is to have a fixed point. The definition of that space has to include the local regularity of the relevant functions $(v, B_2)$ or $(w, G_2)$ and the decay of $(w, G_2)$ as $t\to 0$. As regards the local regularity of $v$ or $w$, we shall use again the space $V$ defined by (3.1), namely we shall take $v$, $w \in {\cal C}(I, V)$ for some interval $I$. On the other hand, when dealing with the Schr\"odinger equation, one time derivative is homogeneous to two space derivatives, so that a given level of regularity in space can be obtained by estimating lower order derivatives if time derivatives are used. In the present case, it turns out that regularity of $v$ is required at the level of $H^k$ with $k > 5/2$, and this is conveniently achieved by using one space and one time derivative, namely with $v \in {\cal C}^1(I, H^1)$. The local regularity of $B_2$ or $G_2$ is less crucial and is essentially dictated by the available estimates. Summing up, the local regularity is conveniently encompassed in the definition of the following space, where $I \subset (0, 1]$~:
$$X_0(I) = \Big \{ (v, B_2) : v \in {\cal C}(I, V) \cap {\cal C}^1(I, H^1 \cap FH^1)\ ,$$
$$B_2 , {\check B}_2 \in {\cal C}(I, \dot{H}^1 \cap \dot{H}^2) \cap {\cal C}^1(I, \dot{H}^1)\Big \} \ . \eqno(4.12)$$

The function space where the map $\Gamma$ is to have a fixed point should also include in its definition the time decay of $(w, G_2)$ as $t \to 0$. For that purpose, we introduce a function $h \in {\cal C}(I, {I \hskip - 1 truemm R}^+)$, where $I = (0, \tau ]$ such that the function $\overline{h}(t) = t^{-3/2} h(t)$ be non decreasing and satisfy
$$\int_0^t dt'\ t{'}^{-1} \ \overline{h}(t') \leq c \ \overline{h}(t)  \eqno(4.13)$$

\noi for some $c > 0$ and for all $t \in I$. Typical examples of such a function would be $h(t) = t^{3/2+ \lambda}(1 - \ell n\ t)^{\mu}$ for some $\lambda > 0$. The function $h(t)$ will characterize the time decay of $\parallel w(t) \parallel_2$ as $t \to 0$. Derivatives of $w$ will have a weaker decay, typically with a loss of one power of $t$ for a time derivative, and of half a power of $t$ for a space derivative. The time decay of $G_2$ is tailored to fit the available estimates, with a loss of one power of $t$ per time or space derivative. Finally, the relevant space for $(w, G_2)$ is defined by 
$$X(I) = \Big \{  (w, G_2) \in X_0 (I) : \parallel (w, G_2);X(I) \parallel \ =\ \mathrel{\mathop {\rm Sup}_{t\in I}}\ h(t)^{-1}\ N (w, G_2)(t) < \infty \Big \}\eqno(4.14)$$ 

\noi where
$$N(w, G_2) = \ \parallel <x> w \parallel_2\ \vee \ t \left ( \parallel <x> \partial_t w \parallel_2 \ \vee\ \parallel <x> \Delta w \parallel_2\right ) $$
$$\vee \ t^{3/2}\left ( \parallel \nabla \partial_t w \parallel_2\ \vee \ \parallel \nabla \Delta w \parallel_2 \right ) \vee \ t^{-1/2} \parallel \nabla G_2 \parallel_2$$
$$\vee \ t^{1/2}\left ( \parallel \nabla^2 G_2  \parallel_2 \ \vee \ \parallel \nabla  \partial_t  G_2  \parallel_2\ \vee \ \parallel \nabla  {\check G}_2 \parallel_2\right )$$
$$\vee \ t^{3/2}\left ( \parallel \nabla^2 {\check G}_2 \parallel_2  \ \vee \ \parallel \nabla \partial_t  G_2 \parallel_2 \right )\ .  \eqno(4.15)$$

We now turn to the implementation of Step 1 and for that purpose we need general boundedness properties of $(v_a, \varphi_a, B_a)$ and general decay assumptions of the remainders $R_i$, $1 \leq i \leq 4$ as $t\to 0$, which we state as follows. Here $I$ is some interval $I = (0, \tau_0]$ with $0 < \tau_0 \leq 1$. \\

\noi (A1) {\bf Boundedness properties of $v_a$.} We assume that 
$$v_a \in ( {\cal C} \cap L^{\infty} ) (I, V) \eqno(4.16)$$
$$t^{1/2} \partial_t v_a \in   ( {\cal C} \cap L^{\infty}) (I, H^2)\ , \ t^{1/2} x\partial_t v_a \in   ( {\cal C} \cap L^{\infty}) (I, H^1) \ . \eqno(4.17)$$
\vskip 5 truemm

\noi (A2) {\bf Boundedness properties of $\varphi_{a}$, $B_a$.} We define $s_a = \nabla \varphi_a$ and we recall that $K_a = B_a + s_a = B_0 + B_{1a} + B_{2a} + s_a$ and that  ${\check B}_{aS} =  {\check B}_{0} +  {\check B}_{1aS} +  {\check B}_{2a}$ in analogy with (2.20 and (2.15). We assume that the following estimates hold for all $t\in I$~:
$$\parallel K_a \parallel_{\infty} \ \leq C(1 - \ell n\ t)\ , \eqno(4.18)$$
$$\parallel \partial_t K_a \parallel_{\infty} \ \vee\ \parallel \nabla K_a \parallel_{\infty} \ \vee \ t \parallel \nabla \partial_t K_a \parallel_{\infty} \ \leq C \ t^{-1}\ , \eqno(4.19)$$
$$\parallel\nabla  s_a \parallel_{\infty} \ \vee\ \parallel \nabla \nabla \cdot s_a \parallel_{3} \ \vee \  \parallel \nabla (B_{1a} + B_{2a}) \parallel_{\infty} $$
$$\vee  \ t \Big ( \parallel \nabla \partial_t s_a \parallel_{\infty}\ \vee \parallel \nabla \partial_t \nabla \cdot s_a \parallel_{3} \ \vee\  \parallel \nabla \partial_t (B_{1a} + B_{2a}) \parallel_{\infty} \Big ) \leq C \ t^{-1/2}\ , \eqno(4.20)$$
$$\parallel \nabla {\check B}_a\parallel_{\infty}\ \vee \ t \parallel\nabla \partial_t {\check B}_a\parallel_{\infty}\ \leq C\ t^{-1}\ , \eqno(4.21)$$
$$\parallel  {\check B}_{aS}\parallel_{\infty}\ \vee \ t \parallel  \partial_t {\check B}_{aS}\parallel_{\infty}\ \leq C\ t^{-1/2}\ . \eqno(4.22)$$ 
\vskip 5 truemm

\noi (A3) {\bf Time decay of the remainders.} We assume that the remainders $R_j$ satisfy the following estimates~:
$$\left ( \parallel <x>R_1(t) \parallel _2\ \leq \right ) \ \parallel <x> \partial_t R_1; L^1((0, t], L^2 ) \parallel \ \leq r_1\ t^{-1}\ h(t)\ , \eqno(4.23)$$
$$\left ( \parallel \nabla R_1(t) \parallel _2\ \leq \right ) \ \parallel \nabla \partial_t R_1; L^1((0, t], L^2 ) \parallel \ \leq r_1\ t^{-3/2}\ h(t)\ ,\eqno(4.24)$$
$$\parallel \nabla^{k+1} R_2(t)  \parallel_2\ \leq r_2\ t^{-1-k\beta}\ h(t)\qquad \hbox{for $k = 0, 1, 2$}\ , \eqno(4.25)$$ 
$$ \parallel \nabla R_3  \parallel_2\ \vee\ t^{1/2}  \parallel \nabla^2 R_3  \parallel_2\ \vee \  t\left ( \parallel \nabla \partial_t  R_3  \parallel_2\ \vee\ \parallel \nabla {\check R}_3  \parallel_2\right )$$
$$\vee\ t^{3/2}  \parallel \nabla^2 {\check R}_3  \parallel_2\ \vee\ t^2  \parallel \nabla \partial_t {\check R}_3  \parallel_2\ \leq r_3 \ h(t)\ ,
\eqno(4.26)$$
$$\parallel \nabla R_4  \parallel_2\ \vee\ t\left (  \parallel \nabla^2 R_4  \parallel_2\ \vee \  \parallel \nabla \partial_t  R_4  \parallel_2\ \vee\ \parallel \nabla {\check R}_4  \parallel_2\right )$$
$$\vee\ t^{2}  \left ( \parallel \nabla^2 {\check R}_4  \parallel_2\ \vee\  \parallel \nabla \partial_t {\check R}_4  \parallel_2\right ) \ \leq r_4\ t^{1/2}\  \ h(t)
\eqno(4.27)$$

\noi for some positive constants $r_j$, $1 \leq j \leq 4$ and for all $t \in I$.\\

The decay properties of the remainders are tailored to fit the definition of the space $X(I)$. For instance when applied to (4.2), the estimate (4.23) yields a contribution $O(h(t))$ to $\parallel <x> w(t) \parallel_2$ by integration over time. \par

In order to ensure that $B_0$ satisfies the properties of $B_a$ appearing in (A2), we need some regularity assumptions on the asymptotic data $(A_+, \dot{A}_+)$. We shall eventually make the following assumption, which suffices for that purpose.\\

\noi (A4) {\bf Regularity of $(A_+, \dot{A}_+)$.} We assume that $\nabla\cdot A_+ = \nabla \cdot \dot{A}_+ = 0$ and that $(A_+, \dot{A}_+)$ satisfies the conditions
$${\cal A}  \in L^2\ , \quad \nabla^2 {\cal A}  \in L^1\ , \quad \omega^{-1} \dot{\cal A} \in L^2\ , \quad \nabla \dot{\cal A} \in L^1 \eqno(4.28)$$
\noi for
$$\left \{ \begin{array}{ll} {\cal A} = (x \cdot \nabla )^{j} \nabla^k A_+ &{\cal A} = (x \cdot \nabla )^{j} \nabla^k (x \cdot A_+)\\ \\ \dot{\cal A}  = (x \cdot \nabla )^{j} \nabla^k \dot{A}_+\qquad\qquad\qquad &\dot{\cal A}  = (x \cdot \nabla )^{j} \nabla^k (x \cdot \dot{A}_+) \end{array} \right . \eqno(4.29)$$ 

\noi for $0 \leq j, k\leq 1$.\\

We can now state the main result of Step 1.\\

\noi {\bf Proposition 4.1.} {\it Let the assumptions (A1)-(A3) be satisfied. Then there exists $\tau$, $0 < \tau \leq 1$, and there exists a unique solution $(w, \psi , G_2)$ of the system (4.2)-(4.5) such that $(w, G_2) \in X((0, \tau ])$ and that $\psi (0) = 0$. Equivalently, there exists a unique solution $(v, \varphi , B_2)$ of the system (2.16)-(2.18) such that $(v - v_a, B_2 - B_{2a}) \in X((0, \tau ])$ and such that $(\varphi - \varphi_a )(t) \to 0$ as $t\to 0$.}\\

We now turn to Step 2, namely the construction of $(v_a, \varphi_a, B_a)$ satisfying the assumptions (A1)-(A3). Taking for orientation $h(t) = t^{3/2 + \lambda}$ with $\lambda > 0$, we need in particular that 
$$\parallel R_1(t) \parallel_2 \ = O \left ( t^{1/2 + \lambda}\right )$$

\noi in order to ensure (4.23). The obvious choice $v_a = U(t) v_+$ yields $\varphi_a = O(1 - \ell n \ t)$ when substituted into (2.17), and therefore $R_1 = O((1 - \ell n\ t)^2)$ when substituted into (4.7), which is off by slightly more than half a power of $t$. We need therefore to go to the next approximation in the resolution of the auxiliary system (2.16)-(2.18). Omitting for the moment the terms containing $B_0$, which require a separate treatment, we choose 
$$\left \{   \begin{array}{ll}
v_a = v_{a0} + v_{a1} \ , &\varphi_a =   \varphi_{a0} + \varphi_{a1}\ , \\
& \\
B_{1a} = B_{1a0} + B_{1a1}\ , &B_{2a} = B_{2a1} + B_{2a2}\ , 
 \end{array} \right . \eqno(4.30)$$

\noi where the last subscript denotes both the order of approximation and the exponent of the power of $t$ as $t \to 0$ (remember that $B_2$ has an extra factor $t$ as compared to $B_1$, cf. (2.11) (2.12)). The lowest approximation is defined by
$$\left \{   \begin{array}{l}
i \partial_t  v_{a0} + (1/2) \Delta v_{a0} = 0 \\
 \\
\partial_t  \varphi_{a0} = t^{-1} g (v_{a0}) + {\check B}_{1L} (v_{a0})\\
 \\
B_{1a0} = B_1 (v_{a0})\\
 \\
B_{2a1} = {\cal B}_2 \left ( v_{a0} , B_{1a0} + \nabla \varphi_{a0}\right )
 \end{array} \right .\eqno(4.31)$$ 
 
 \noi with initial conditions $v_{a0} (0) = v_+$, $\varphi_{a0} (1) = 0$, so that in particular $v_{a0} = U(t) v_+$ and $\varphi_{a0} = O(1 - \ell n\ t)$. Similarly the second approximation is defined by 
$$\left \{   \begin{array}{l}
i \partial_t  v_{a1}  = \hbox{remaining zero order terms of (2.16) not containing $B_0$} \\
 \\
\partial_t  \varphi_{a1} = 2t^{-1} g (v_{a0}, v_{a1}) + 2{\check B}_{1L} (v_{a0}, v_{a1})\\
 \\
B_{1a1} = 2 B_1 (v_{a0}, v_{a1})\\
 \\
B_{2a2} = \hbox{order one terms of ${\cal B}_2 \left ( v_{a} , B_{1a} + B_{2a} \nabla  \varphi_{a}\right )$}
 \end{array} \right .\eqno(4.32)$$ 

\noi with initial conditions $v_{a1} (0) = 0$, $\varphi_{a1}(0) = 0$. It turns out that the previous choice yields $v_{a1} = O(t(1 - \ell n \ t)^2)$, $\varphi_{a1} = O(t(1 - \ell n\ t)^2)$, and suffices to control the $B_0$ independent part of the remainders.\par

We now turn to the terms containing $B_0$, which are not covered by the previous choice. The difficulty comes from the fact that applying a derivative to $B_0$, whether in space or in time, generates a factor $t^{-1}$, as can be seen on the change of variable (2.6). Thus the most dangerous terms come from $\nabla \partial_t R_1$ which contains typically $(\nabla \partial_t B_0)\cdot \nabla v_+$. Taking again $h(t) = t^{3/2+ \lambda }$ for orientation, one needs an estimate
$$\parallel \nabla \partial_t R_1\parallel_2\  \leq O\left ( t^{-1 + \lambda } \right )$$

\noi whereas the best estimates on $B_0$ for $A_0$ a solution of the free wave equation yield only 
$$\parallel ( \nabla \partial_t B_0)\cdot \nabla v_+\parallel_2\ \leq \ \parallel \nabla \partial_t B_0\parallel_2\ \parallel \nabla v_+ \parallel_{\infty} \ \leq O\left ( t^{-3/2} \right )$$

\noi which is again off by slightly more than half a power of $t$. In order to deal with that problem, we impose a support condition on $v_+$, namely 
$${\rm Supp}\ v_+ \subset \left \{ x: |\ |x|-1| \geq \eta \right \} \eqno(4.33)$$

\noi for some $\eta > 0$. The effect of that condition is best understood when combined with $B_0$ coming from a solution $A_0$ of the free wave equation with compactly supported data. In fact, if 
$${\rm Supp}\ (A_+, \dot{A}_+) \subset \left \{ x : |x| \leq R \right \}\ ,$$

\noi then by the Huyghens principle, $A_0$ defined by (2.2) satisfies
$${\rm Supp}\ A_0 \subset \left \{ (x, t) : |\ |x|-t| \leq R \right \}$$

\noi so that by (2.6)$_0$
$${\rm Supp}\ B_0 \subset \left \{ (x, t) : |\ |x|-1| \leq t\ R \right \}$$

\noi and therefore for all integers $j$, $k$, $\ell$, 
$$\left ( \nabla^k \partial_t^j B_0\right ) \nabla^{\ell} v_+ = 0 \quad \hbox{for $t < \eta /R$}\ .$$

\noi More generally, the support condition (4.33) together with some decay of $(A_+, \dot{A}_+)$ for large $|x|$ suffices to ensure the required estimates. We shall make the following assumption.\\\

\noi (A5) {\bf Space decay of $(A_+, \dot{A}_+)$.} Let $\chi_R$ be the characteristic function of the set $\{x:|x|\geq R$). We assume that $(A_+, \dot{A}_+)$ satisfies
$$\left \{ \begin{array}{l} \parallel \chi_R \nabla^k (x\cdot \nabla )^j A_+ \parallel_2\ \vee\ \parallel \chi_R \nabla^k (x\cdot \nabla)^j x\cdot A_+ \parallel_2\ \leq\ C\ R^{- 1}\\ \\
\parallel \chi_R (x\cdot \nabla )^j \dot{A}_+;L^2 \cap L^{6/5} \parallel\ \vee\ \parallel \chi_R (x \cdot  \nabla )^j x\cdot \dot{A}_+;L^2 \cap L^{6/5}  \parallel\ \leq\ C\ R^{- 1} \end{array}\right . \eqno(4.34)$$

\noi for $0 \leq j , k \leq 1$ and for all $R \geq R_0$ for some $R_0 > 0$. \\

The support condition (4.33) appeared in the early work \cite{11r} on the MS system and was subsequently eliminated \cite{2r} \cite{10r} in the framework of the method of \cite{11r}. It is an open question whether that condition can also be eliminated in the framework of the more complicated method described in the present paper.\par

We can now state the result of the completion of Step 2. \\

\noi {\bf Proposition 4.2.} {\it Let $v_+ \in H^5$, $xv_+ \in H^4$ and let $v_+$ satisfy the support condition (4.33). Let $(A_+, \dot{A}_+)$ satisfy the regularity and decay assumptions (A4) (A5). Let $(v_a, \varphi_a , B_a)$ be defined by (4.30)-(4.32) with $v_{a0}(0) = v_+$, $\varphi_{a0}(1) = 0$, $v_{a1}(0) = 0$, $\varphi_{a1}(0) = 0$ and $B_{0a} = B_0$. Then the assumptions (A1)-(A3) are satisfied with} 
$$h(t) = t^2 (1 - \ell n\ t)^4 \ . \eqno(4.35)$$

Putting together Propositions 4.1 and 4.2 yields the main result on the Cauchy problem at $t = 0$ for the auxiliary system (2.16)-(2.18).\\

\noi {\bf Proposition 4.3.} {\it Let $v_+ \in H^5$, $xv_+ \in H^4$ and let $v_+$ satisfy the support condition (4.33). Let $(A_+, \dot{A}_+)$ satisfy the regularity and decay assumptions (A4) (A5). Let $(v_a, \varphi_a , B_a)$ be defined by (4.30)-(4.32) with $v_{a0}(0) = v_+$, $\varphi_{a0}(1) = 0$, $v_{a1}(0) = 0$, $\varphi_{a1}(0) = 0$ and $B_{0a} = B_0$. Let $h(t) = t^2 (1 - \ell n\ t)^4$. Then there exists $\tau$, $0 < \tau \leq 1$, and there exists a unique solution $(v, \varphi, B_2)$ of the auxiliary system (2.16)-(2.18) such that $(v-v_a, B_2- B_{2a}) \in X((0, \tau ])$ and such that $(\varphi - \varphi_a)(t) \to 0$ as $t \to 0$.}\\

We now return to the Cauchy problem at infinity for the MS system in the form (1.3) (2.1). We start from the asymptotic data $(u_+, A_+, \dot{A}_+)$ for $(u, A)$ and we define $v_+ = \overline{Fu_+}$. We define $A_0$ by (2.2) and $B_0$ by (2.6)$_0$. We define $(v_a, \varphi_a, B_a)$ by (4.30)-(4.32) with the appropriate initial conditions. We then define the asymptotic form $(u_a, A_a)$ of $(u, A)$ in analogy with (2.4)-(2.6), (2.14) by 
$$u_{ac} = v_a \exp \left ( - i \varphi_a\right )\ , \eqno(4.36)$$
$$\widetilde{u}_a(t) = \overline{F \widetilde{u}_{ac} (1/t)}\ , \eqno(4.37)$$
$$A_a(t) = - t^{-1} D_0 (t)\ B_a (1/t) \ . \eqno(4.38)$$

The final result can then be stated as follows. We recall that $\widetilde{u}$ is defined by (2.5), that $V_{\star}$ is defined by (3.2), that $S$ is defined by (3.3) and that $j$, $k$, $\ell$ denote non negative integers.\\
 
\noi {\bf Proposition 4.4.} {\it Let $u_+$ be such that $v_+ \equiv \overline{Fu_+} \in H^5$ with  $xv_+\in H^4$ and that $v_+$ satisfy the support condition (4.33). Let $(A_+, \dot{A}_+)$ satisfy the regularity and decay assumptions (A4) (A5). Let $(u_a, A_a)$ be defined by (4.36)-(4.38) with $(v_a, \varphi_a, B_a)$ defined by (4.30)-(4.32) with $v_{a0}(0) = v_+$, $\varphi_{a0}(1) = 0$, $v_{a1}(0) = 0$, $\varphi_{a1}(0) = 0$ and $B_{0a} = B_0$. Then there exists $T \geq 1$ and there exists a unique solution $(u, A)$ of the MS system (1.3) (2.1) in $I = [T , \infty )$ with the following properties~: 
$$\widetilde{u} \in {\cal C}(I, V_{\star})\quad , \quad \partial_t \widetilde{u} \in {\cal C} (I, H^1 \cap F H^1)\ ,$$
$$ A \in {\cal C} (I , \dot{H}^1 \cap \dot{H}^2) \ , \ x\cdot A \in {\cal C} (I , \dot{H}^2)\ , \ SA \in {\cal C} (I, H^1)\ ,$$

\noi and $(u, A)$ satisfy the estimates
$$\parallel \widetilde{u}(t) ; V_{\star} \parallel \ \leq C (1 + \ell n\ t)^3 \eqno(4.39)$$
$$\parallel  \nabla (A - A_0)(t) \parallel_2 \ \leq C\ t^{-1/2} \eqno(4.40)$$

\noi for all $t\in I$. Furthermore $(u,A)$ behaves asymptotically as $(u_a, A_a)$ in the sense that the following estimates hold for all $t\in I$~:
$$\parallel x^k \nabla^{\ell} ( \widetilde{u} -  \widetilde{u}_a)(t)\parallel_2\ \leq C \ t^{-2 + k/2} (1 + \ell n\ t)^4 \eqno(4.41)$$

\noi for $\ell = 0, 1$ and $0 \leq k + \ell \leq 3$, 
$$\parallel x^k \partial_t ( \widetilde{u} -  \widetilde{u}_a)(t)\parallel_2\ \leq C \ t^{-3 + k/2} (1 + \ell n\ t)^4 \eqno(4.42)$$

\noi for $k = 0, 1$, 
$$\parallel \nabla^{k+1} S^j (A - A_a) (t)\parallel_2 \ \vee \ t^{-1} \parallel \nabla^{k+1} S^j x \cdot (A - A_a) (t) \parallel_2$$
$$\leq C\ t^{-5/2-k/2} (1 + \ell n\ t)^4$$
\noi for $0 \leq j, k, j+k \leq 1$.}\\

The existence part of the proof of Proposition 4.4 is obtained by first constructing $(v, \varphi , B_2)$ by the use of Proposition 4.3, then constructing $(u, A)$ from $(v, \varphi , B_2)$ by using the change of variables (2.4)-(2.6) (2.10) (2.11) (2.14), and finally translating the properties of $(v, \varphi , B_2)$ that follow from Proposition 4.3 in terms of $(u, A)$. The uniqueness part of the proof follows readily from Proposition 3.2. In fact, if $(u_i, A_i)$, $i = 1,2$, are two solutions of the system (1.3) (2.1) obtained from Proposition 4.4 with the same $(u_a, A_a)$, then the conditions (3.9) (3.10) follow from (4.39) (4.40) with $\alpha = 3$, while the condition (4.41) together with the commutation relation $xU(-t) = U(-t) (x + it \nabla )$ implies 
$$\parallel  <x/t> (u_1-u_2)(t) \parallel_2 \ \leq C\ t^{-2} (1 + \ell n\ t)^4$$

\noi which implies (3.11).

\end{document}